\documentclass[12pt,a4paper]{article}
\usepackage{amsfonts}

\usepackage{amssymb}

\usepackage{amsmath}

\newtheorem{theorem}{Theorem}[section]

\newtheorem{lemma}[theorem]{Lemma}

\newtheorem{proposition}[theorem]{Proposition}

\date{}

\begin{document}

\title{ Composition-Diamond lemma for $\lambda$-differential  associative algebras with
multiple operators\footnote{Supported by the NNSF of China
(No.10771077) and the NSF of Guangdong Province (No.06025062).} }

\author{Jianjun Qiu\\
{\small \ Mathematics and Computational  Science School }\\
{\small \ Zhanjiang  Normal University}\\
{\small \ Zhanjiang  524048, China}\\
{\small \ jianjunqiu@126.com}\\
\\
Yuqun Chen\\
{\small \ School of Mathematical Sciences, South China Normal
University}\\
{\small Guangzhou 510631,   China}\\
{\small  yqchen@scnu.edu.cn}\\
\\
} \vspace{4mm}

\maketitle \noindent\textbf{Abstract:} In this paper, we establish
the Composition-Diamond lemma for $\lambda$-differential associative
algebras over a field $K $ with multiple operators. As applications,
we obtain Gr\"{o}bner-Shirshov bases of free $\lambda$-differential
Rota-Baxter algebras.  In particular, linear bases of free
$\lambda$-differential Rota-Baxter algebras are obtained and
consequently, the free $\lambda$-differential Rota-Baxter algebras
are constructed by words.



\section{Introduction}
 Let $K$ be a field  and $\lambda\in K$. A differential
algebra of weight $\lambda$, namely,  a $\lambda$-differential
algebra (see \cite{lw, ko}), is an associative $K$-algebra $R$ with
a $\lambda$-differential operator $D:R\rightarrow R $ such that
$$
D(xy)=D(x)y+xD(y)+\lambda D(x)D(y), \forall x, y\in R.
$$
A Rota-Baxter algebra of weight $\lambda$ (see \cite{Bax, ro}) is an
associative $K$-algebra $R$ with a Rota-Baxter operator
$P:R\rightarrow R$ such that
$$
P(x)P(y)=P(xP(y))+P(P(x)y)+ \lambda P(xy),\forall x, y\in R.
$$
Similar to the relation of integral and differential operators, L.
Guo and W. Keigher \cite{lw} introduced the notion of
$\lambda$-differential Rota-Baxter algebra which is a $K$-algebra
$R$ with  a $\lambda$-differential operator $D$ and a Rota-Baxter
operator $P$ such that $DP=Id_{R}$.

There have been  some constructions of free Rota-Baxter algebras
(commutative and associative). We note that G.-C. Rota \cite{ro} and
P. Cartier \cite{ca} gave the explicit constructions of the free
commutative Rota-Baxter algebras on a set when $\lambda=1$, namely,
the shuffle Baxter and standard Baxter algebras, respectively.
Recently, L. Guo  and W. Keigher \cite{gk, gk1} constructed the free
commutative Rota-Baxter algebras (with identity or without identity)
for any $\lambda \in K$ by the mixable shuffle product. These
algebras are now called the mixable shuffle product algebras. In
fact, these algebras generalize the classical construction of
shuffle product algebras. K. Ebrahimi-Fard and L. Guo \cite{kl1}
constructed recently the free associative Rota-Baxter algebras  by
Rota-Baxter words.

E.  Kolchin \cite{ko} considered  the differential algebra and
constructed a free   differential algebra. L. Guo and W. Keigher
\cite{lw} dealt with a generalization  of this algebra. Also in
\cite{lw}, the free $\lambda$-differential Rota-Baxter algebra was
obtained by using the free Rota-Baxter algebra  on  planar decorated
rooted trees.

K. Ebrahimi-Fard and L. Guo  \cite{kl2} used rooted trees and
forests to give an explicit construction of free noncommutative
Rota-Baxter algebras on modules and sets. K. Ebrahimi-Fard, J. M.
Gracia-Bondia and F. Patras \cite{EGP} gave the  solution of the
word problem for free non-commutative Rota--Baxter algebra. A free
Rota--Baxter algebra was constructed on decorated trees by M. Aguiar
and M. Moreira \cite{AM}.

The concept of multi-operators algebras ($\Omega$-algebras) was
first introduced by A. G. Kurosh  in \cite{ku, ku2}. Also, Kurosh
noticed that free $\Omega$-algebras share many of the combinatorial
properties of free non-associative algebras. On the other hand, the
Gr\"{o}bner-Shirshov bases theory for Lie algebras was first
considered by A. I. Shirshov \cite{S3}. In fact, Shirshov \cite{S3}
defined the composition of two Lie polynomials and established the
Composition lemma for the Lie algebras.  Later on, L. A. Bokut
\cite{b76} specialized the approach of Shirshov to associative
algebras, see also G. M. Bergman \cite{b}. For commutative
polynomials, this lemma is known as the Buchberger's Theorem in
\cite{bu65, bu70}.

Gr\"{o}bner-Shirshov bases for $\Omega$-algebras were given in the
paper of V. Drensky and R. Holtkamp \cite{dr}. In a recent paper of
L. B. Bokut, Y. Chen and J. Qiu \cite{b08}, the Composition-Diamond
lemma is established for associative $\Omega$-algebras.

In this paper, we construct free $\lambda$-differential associative
algebras with multiple operators and give the Composition-Diamond
lemma for such algebras.  As applications, we obtain
Gr\"{o}bner-Shirshov bases of free $\lambda$-differential
Rota-Baxter algebras. Then, linear bases of free
$\lambda$-differential Rota-Baxter algebras are abtained, which is
the same that obtained by Bokut,  Chen and Qiu \cite{b08}, and
similar to those obtained by using other methods in the paper of Guo
and Keigher \cite{lw}.

\section{ Free  $\lambda$-differential associative algebras with
multiple operators}

Let $K$ be a field and $\lambda \in K$. A $\lambda$-differential
associative algebra with multiple operators  is a
$\lambda$-differential algebra $R$ with a set
 $\Omega$ of multi-linear  operators. In order to construct the
free $\lambda$-differential associative algebras with multiple
operators,  we recall the construction of free
$\lambda$-differential algebras firstly. For the detail of the free
$\lambda$-differential algebras, see \cite{lw, b08,ko}.

Let $X$ be a set, $D$ a symbol and
$$
D^{\omega}(X)=\{D^i(x)|i\geq 0, x\in X\}
$$
 where $D^0(x)=x$.

Let $S(D^{\omega}(X))$ be  the free  semigroup generated by
$D^{\omega}(X)$. Denote $KS(D^{\omega}(X))$ the semigroup algebra
over the semigroup $S(D^{\omega}(X))$.  Extend linearly $D: \
KS(D^{\omega}(X))\rightarrow KS(D^{\omega}(X))$ in  the following
way: for any $u=u_1u_2\cdots u_t \in S(D^{\omega}(X))$, where
$u_i\in D^{\omega}(X)$, define $D(u)$ by induction on $t$:

If $t=1$, i.e., $u=D^i(x)$ for some $i\geq 0, x\in X $, then
$$
D(u)=D^{i+1}(x).
$$

If $t>1$, then
$$
D(u)=\lambda D(u_1)D(u_2\cdots u_t)+D(u_1)(u_2\cdots u_t
)+u_1D(u_2\cdots u_t).
$$
\begin{theorem}\label{t1}{\em (\cite{ b08,lw,ko})}
$(KS(D^{\omega}(X)), D)$ is a free $\lambda$-differential algebra on
$X$.\hfill $\blacksquare$
\end{theorem}

Secondly, we construct free $\lambda$-differential associative
algebras with multiple operators.

Let
$$
\Omega=\bigcup_{n=1}^{\infty}\Omega_{n}
 $$
where $\Omega_{n}$ is the set of $n$-ary operators,  ary
$(\delta)=n$ if $\delta\in \Omega_n$.

Define
$$\mathfrak{L}_{0}=S(D^{\omega}(X_0)), \ \ X_0=X$$
$$\mathfrak{L}_{1}=S(D^{\omega}(X_{1})), \ \ X_1=X\cup \Omega(\mathfrak{L}_{0})
$$
where
$$
\Omega(\mathfrak{L}_{0})=\bigcup\limits_{t=1}^{\infty}\{\delta(y_1,y_2,\cdots,y_t)|\delta\in
\Omega_t, y_i\in\mathfrak{L}_{0}, \ i=1,2,\cdots,t\}.
$$
For $n>1$, define
$$
\mathfrak{L}_{n}=S(D^{\omega}(X_{n})), \ \ X_n=X\cup
\Omega(\mathfrak{L}_{n-1})
$$
where
$$
\Omega(\mathfrak{L}_{n-1})=\bigcup\limits_{t=1}^{\infty}\{\delta(y_1,y_2,\cdots,y_t)|\delta\in
\Omega_t, y_i\in \mathfrak{L}_{n-1}, \ i=1,2,\cdots,t\}.
$$
Then we have
$$
\mathfrak{L}_{0}\subset\mathfrak{L}_{1}\subset\cdots
\subset\mathfrak{L}_{n}\subset\cdots.
$$
Let
$$
\mathfrak{L}(X)=\bigcup_{n\geq0}\mathfrak{L}_{n}.
$$
Then, it is easy to see that $\mathfrak{L}(X)$ is a semigroup such
that $ \Omega(\mathfrak{L}(X))\subseteq \mathfrak{L}(X)$.

 For any
$v\in \mathfrak{L}(X)$, $v$ has  a unique expression
$$
v=v_1v_2\cdots v_n
$$
where each $v_i\in D^{\omega}(X\cup \Omega(\mathfrak{L}(X))).$  If
this is the case, then we define $bre(v)=n$.

Let  $D\langle X;\Omega\rangle=K\mathfrak{L}(X)$ be the semigroup
algebra over $\mathfrak{L}(X)$. Then, similar to Theorem \ref{t1},
$D\langle X;\Omega\rangle$ is a $\lambda$-differential algebra.

For any $\theta\in\Omega_n$, define
$$\theta:\ \mathfrak{L}(X)^n\rightarrow \mathfrak{L}(X), \
(v_1,v_2,\cdots,v_n)\mapsto \theta(v_1,v_2,\cdots,v_n)
$$
and extend linearly  $\theta: \ D\langle
X;\Omega\rangle^n\rightarrow D\langle X;\Omega\rangle$. Then it is
easy to see that $D\langle X;\Omega\rangle$ is a free
$\lambda$-differential associative algebra on $X$ with multiple
operators $\Omega$.

Now, we describe the elements in $ \mathfrak{L}(X)$  by  labeled
reduced planar rooted forests, see \cite{dr, kl2}.

For any $x\in X_0$,  $x$ can  be described by a labeled reduced
planar rooted tree $\bullet_{x}$. Then for any  element
$v=D^{i}(x),\ i\geq 1$ can be described by

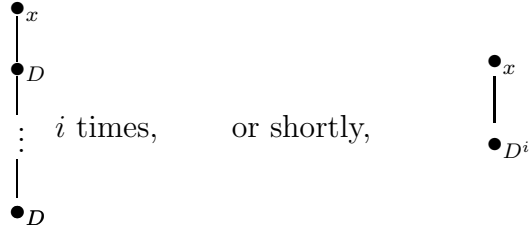
\setlength{\unitlength}{1cm}
\begin{picture}(4,3.0)
 \put(4,0.4){\line(0,1){0.5}}
\put(4,1.5){\line(0,1){0.5}} \put(4,2.2){\line(0,1){0.6}}\put(3.9,
0.1){$\bullet_{D}$} \put(4,1.0){$\vdots$} \put(3.9,
0.1){$\bullet_{D}$}\put(3.9, 2.0){$\bullet_{D}$}\put(3.9,
2.8){$\bullet_{x}$}\put(4.5 ,1.2){$ i \mbox{ times, \ \ \ \ \ \ or
shortly, }$}
\end{picture}
\setlength{\unitlength}{1cm}
\begin{picture}(4,3.0)\put(6,1.4){\line(0,1){0.6}}
\put(5.9, 1){$\bullet_{D^i}$}\put(5.9, 2.1){$\bullet_{x}$}
\end{picture}\\

Then for  $u=u_1 u_2\dots u_n\in S(D^{\omega}(X))$, we  describe $u$
 by
$$
\bullet_{u_1}\sqcup \bullet_{u_2}\sqcup \dots \sqcup\bullet_{u_n}
$$
which is called a  labeled reduced planar rooted forest. For
example, if $u=D^{3}(x_1)D^2(x_2)$, then  $u$ can be described by

 \setlength{\unitlength}{1cm}
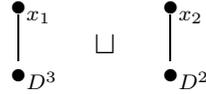
\begin{picture}(4,1.2)
\put(5,0.4){\line(0,1){0.6}}\put(4.9, 0.1){$\bullet_{D^3}$}\put(4.9,
1){$\bullet_{x_1}$}\put(6.0, 0.5){$\sqcup$}
\put(7,0.4){\line(0,1){0.6}}\put(6.9, 0.1){$\bullet_{D^2}$}
 \put(6.9, 1){$\bullet_{x_2}$}
\end{picture}\\
For any $u\in \Omega(\mathfrak{L}_{0})$, say $u=\theta_n(u_1, \dots
,u_n)$ where each $u_i\in \mathfrak{L}_{0}$, $u$ can be described by

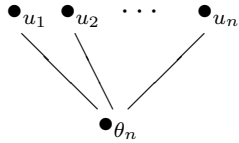
\setlength{\unitlength}{1cm}
\begin{picture}(4,2.0)
\put(7,0.5){\line(1,1){1}}
\put(6.8,0.5){\line(-1,2){0.5}}
\put(6.6,0.5){\line(-1,1){1}}
\put(6.6,0.2){$\bullet_{\theta_n}$}
\put(6.1,1.7){$\bullet_{u_2}$}
\put(5.4,1.7){$\bullet_{u_1}$}
\put(7.9,1.7){$\bullet_{u_n}$}\
\put(6.9,1.7){$\cdots$}
\end{picture}\\
For example,  if $u=\theta_3(u_1, u_2,u_3)\in
\Omega(\mathfrak{S}_{0})$ where $u_1=D(x_1) D^2(x_2), \
u_2=D^2(x_3), \ u_3=D^3(x_4)D^3(x_5)$, then
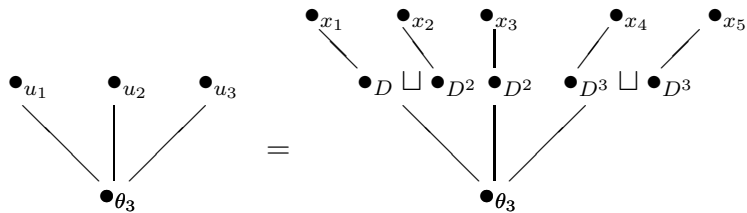
\begin{eqnarray*}
\setlength{\unitlength}{1cm}
\begin{picture}(4,3.0)
 \put(1,0.5){\line(1,1){1}}
 \put(0.8,0.5){\line(0,1){1}}
\put(0.6,0.5){\line(-1,1){1}} \put(0.6,0.2){$\bullet_{\theta_3}$}
\put(0.6,0.2){$\bullet_{\theta_3}$} \put(0.7,1.7){$\bullet_{u_2}$}
\put(-0.6,1.7){$\bullet_{u_1}$} \put(1.9,1.7){$\bullet_{u_3}$}
\put(2.8, 0.8){=}
\end{picture}
\setlength{\unitlength}{1cm}
\begin{picture}(4,3.0)
\put(2,0.5){\line(1,1){1}} \put(1.8,0.5){\line(0,1){1}}
\put(1.6,0.5){\line(-1,1){1}} \put(1.6,0.2){$\bullet_{\theta_3}$}
\put(1.6,0.2){$\bullet_{\theta_3}$} \put(1.7,1.7){$\bullet_{D^2}$}
\put(0,1.7){$\bullet_{D}\sqcup\bullet_{D^2}$}
\put(2.7,1.7){$\bullet_{D^3}\sqcup\bullet_{D^3}$}
\put(0.0,2.0){\line(-1,1){0.5}} \put(-0.7,2.6){$\bullet_{x_1}$}
\put(1,2.0){\line(-3,4){0.4}} \put(0.5,2.6){$\bullet_{x_2}$}
\put(1.8,2.0){\line(0,1){0.5}} \put(1.6,2.6){$\bullet_{x_3}$}
\put(2.9,2.0){\line(3,4){0.4}} \put(3.3,2.6){$\bullet_{x_4}$}
\put(4.0,2.0){\line(1,1){0.5}} \put(4.6,2.6){$\bullet_{x_5}$}
\end{picture}
\end{eqnarray*}

 Now, for any  $v=v_1v_2\cdots v_n\in
\mathfrak{L}_{1}=S(D^{\omega}(X\cup \Omega(\mathfrak{L}_{0}))) \
(v_i\in D^{\omega}(X\cup \Omega(\mathfrak{L}_{0})), \ i=1,\dots,n)$,
$v$ can be described by labeled reduced planar rooted forest:
$$
\bullet_{v_1}\sqcup\bullet_{v_2}\sqcup\cdots\sqcup\bullet_{v_n}.
$$
With the above way, we can describe all the elements in
$\mathfrak{S}(X)$ by labeled reduced planar rooted forests.

Therefore,  each element in $\mathfrak{S}(X)$ corresponds uniquely
to a labeled reduced planar rooted forest.
\\

For $1\leq t\leq n$, denote by
$$
\Theta_{n}^{t}=\{(i_1,\cdots,i_n)\in \{0,
1\}^{n}|i_1+\cdots+i_n=t\}.
$$

The following proposition is useful in the next section.

\begin{proposition}\label{p2.2}
Let  $v_t\in D^{\omega}(X\cup \Omega(\mathfrak{L}(X))), t=1,\cdots,
n$. Then
\begin{eqnarray*}
D(v_1v_2\cdots v_n)&=&\sum\limits_{(i_1,\cdots,
i_n)\in\Theta_{n}^{1}}D^{i_1}(v_1)\cdots
D^{i_n}(v_n)\\
&&+\sum_{t=2}^{n} \sum\limits_{(i_1,\cdots,
i_n)\in\Theta_{n}^{t}}\lambda ^{t-1}D^{i_1}(v_1)\cdots D^{i_n}(v_n).
\end{eqnarray*}
\end{proposition}
{\bf Proof:} If $\lambda=0$, the result is clear. Assume that
$\lambda\neq0$. Induction on $n$. For $n=1,2$, the result is clear.
Now we assume that for $n-1$ the result  is true. Since
$$
D(v_1v_2\cdots v_n)=\lambda D(v_1)D(v_2\cdots v_n)+D(v_1)v_2\cdots
v_n +v_1D(v_2\cdots v_n)
$$
and by induction, we have
\begin{eqnarray*}
D(v_1v_2\cdots v_n)
&=&\sum_{t=2}^{n} \sum\limits_{(1,i_2,\cdots,
i_n)\in\Theta_{n}^{t}}\lambda ^{t-1} D(v_1)D^{i_2}(v_2)\cdots
D^{i_n}(v_n)\\
&&+ D(v_1)D^{0}(v_2)\cdots D^{0}( v_n)\\
&&+\sum_{t=1}^{n-1} \sum\limits_{(0,i_2,\cdots,
i_n)\in\Theta_{n}^{t}}\lambda ^{t-1}D^{0}(v_1)D^{i_2}(v_2)\cdots
D^{i_n}(v_n)\\
&=& \lambda^{n-1} D(v_1)\cdots D(v_n)\\
&&+\sum_{t=1}^{n-1}( \sum\limits_{(1,i_2,\cdots,
i_n)\in\Theta_{n}^{t}}\lambda ^{t-1} D(v_1)D^{i_2}(v_2)\cdots
D^{i_n}(v_n)\\
&&+ \sum\limits_{(0,i_2,\cdots, i_n)\in\Theta_{n}^{t}}\lambda
^{t-1}D^{0}(v_1)D^{i_2}(v_2)\cdots
D^{i_n}(v_n))\\
&=&\sum_{t=1}^{n} \sum\limits_{(i_1,\cdots,
i_n)\in\Theta_{n}^{t}}\lambda ^{t-1}D^{i_1}(v_1)\cdots D^{i_n}(v_n).
\ \ \ \ \ \ \ \ \ \ \  \ \ \ \ \ \ \ \ \ \ \ \ \ \  \blacksquare
\end{eqnarray*}

\section{ Composition-Diamond lemma for $\lambda$-differential  associative algebras with
multiple operators}

In this section, we introduce the notions of Gr\"{o}bner-Shirshov
bases  for $\lambda$-differential  associative algebras with
multiple operators and establish the Composition-Diamond lemma for
such algebras.

Let $X$ be a set  and $D\langle X;\Omega\rangle$ the  free
$\lambda$-differential  associative algebra with multiple operators
$\Omega$.  For any $u\in \mathfrak{L}(X)$,
$$dep(u)=\mbox{min}\{n|u\in\mathfrak{L}_{n} \}
$$
is called the depth of  $u$. Let $deg_{_{\Omega}}(u)$ be the number
of the occurrences of $\theta \in \Omega$ in $u$, for example, if
$u=\theta_t(u_1,\cdots,u_t)\in \Omega(\mathfrak{L}(X))$, then
$deg_{_{\Omega}}(u)=1+\sum_{i=1}^tdeg_{_{\Omega}}(u_i)$.

Let $X$ and $\Omega$ be well ordered. We define an ordering $>$ on
$\mathfrak{L}(X)=\bigcup_{n\geq0}\mathfrak{L}_{n}$ by induction on
$n$.

Suppose that $n=0$, $\mathfrak{L}_{0}=S(D^{\omega}(X_0))$.

Order $D^{\omega}(X_{0})$ by
\begin{equation*}\label{e1}
D^{i}(u)>D^{j}(v)\ \mbox{ if }\ (i,u)>(j,v)\ \
 \mbox{ lexicographically}.
\end{equation*}

Order $\mathfrak{L}_{0}=S(D^{\omega}(X_{0}))$. For any $u=u_1\cdots
u_t$,
 $v=v_1\cdots v_s\in \mathfrak{L}_{0}$, (i.e., $u_i,v_j\in D^{\omega}(X_0)$), define
 $u> v$ if

\begin{equation*}\label{e2}
(bre(u), u_1,\cdots, u_t)>
 (bre(v), v_1,\cdots, v_s)\ \
 \mbox{ lexicographically}.
\end{equation*}

Suppose that $n>0$. We order $\mathfrak{L}_{n}$ in three steps.

Firstly, order $X_n=X\cup
\Omega(\mathfrak{L}_{n-1})$ by \\
(i) \ $x<u$ if $u\in \Omega(\mathfrak{L}_{n-1})$ and $x\in X$; \\
(ii) \ for $u=\theta_t(u_1,\cdots,u_t),v=\delta_s(v_1,\cdots,v_s)\in
\Omega(\mathfrak{L}_{n-1})$,  $u>v$ if
$$
(deg_{_{\Omega}}(u),\theta_t, u_1 \cdots,
u_t)>(deg_{_{\Omega}}(v),\delta_s, v_1, \cdots,
 v_s)\ \
 \mbox{ lexicographically}.
$$

Secondly, order $D^{\omega}(X_n)$. For any $u, v\in X_n$, define
\begin{equation*}\label{e1}
D^{i}(u)>D^{j}(v)\ \mbox{ if }\
(deg_{_{\Omega}}(u),i,u)>(deg_{_{\Omega}}(u),j,v)\ \
 \mbox{ lexicographically}.
\end{equation*}

Thirdly, order $S(D^{\omega}(X_n))$. For any $u=u_1\cdots u_t,\
 v=v_1\cdots v_s\in S(D^{\omega}(X_n))$, define  $u> v $  if
\begin{equation*}\label{e5}
(deg_{_{\Omega}}(u),bre(u), u_1,\cdots, u_t)>
 (deg_{_{\Omega}}(v),bre(v), v_1,\cdots, v_s)\ \
 \mbox{ lexicographically}.
\end{equation*}

It is easy to check that $>$ is a well ordering on
$\mathfrak{L}(X)$. Throughout this paper,  this ordering will be
used.

Then for any $0\neq f\in D\langle X;\Omega\rangle$, $f$ has a
leading term $\bar{f}$ and
$$
f=\alpha_1\bar{f} +\sum_{i=2}^{m}\alpha_iu_i
$$
where $\bar{f}>u_i$, $i=2,\dots, m$ and $0\neq \alpha_i \in K$.
Denote by  $lc(f)$  the coefficient of the leading term $\bar{f}$.
If $lc(f)$ equals $1$, we call $f$  monic.

The proof of the following lemma is straightforward.

\begin{lemma}\label{l3.1}
Let $u,v\in \mathfrak{L}(X),u>v, \ \theta_n\in \Omega_n$ and $w_j\in
\mathfrak{L}(X), \ j=1,\dots,n$. Then,
\begin{enumerate}
\item[(i)]
for any $i, \ 1\leq i\leq n$,
$$
\theta_n(w_1,\dots,w_{i-1},u,w_{i+1},\dots,
w_n)>\theta_n(w_1,\dots,w_{i-1},v,w_{i+1},\dots, w_n);
$$
\item[(ii)] for any $a,b \in \mathfrak{L}(X)$,
$$
au>av, \ ub>vb.
$$\hfill $\blacksquare$
\end{enumerate}
\end{lemma}

\begin{lemma} \label{l3.3}
Let  $v_t\in D^{\omega}(X\cup \Omega(\mathfrak{L}(X))),t=1,\cdots,
n$.

(i) If  $\lambda\neq 0$, then
$$
\overline{D^{i}(v_1v_2\cdots v_n)}=D^{i}(v_1)D^{i}(v_2)\cdots
D^{i}(v_n)
$$
and the coefficient of $\overline{D^{i}(v_1v_2\cdots v_n)}$ is
$\lambda^{(n-1)i}$.

(ii)  If  $\lambda= 0$, then
$$
\overline{D^{i}(v_1v_2\cdots v_n)}=D^{i}(v_1)v_2\cdots v_n
$$
and the coefficient of $\overline{D^{i}(v_1v_2\cdots v_n)}$ is 1.

It follows that if $u,v \in \mathfrak{L}(X)$ and $u>v$, then
$\overline{D(u)}>\overline{D(v)}$.
\end{lemma}
{\bf Proof:} Clearly,  by using induction on $i$, (ii) follows from
Proposition \ref{p2.2}. Now, we prove (i).

In Proposition \ref{p2.2}, we have
\begin{equation}\label{e3}
D(v_1v_2\cdots v_n)=\lambda^{n-1}D(v_1)\cdots
D(v_n)+\sum_{t=1}^{n-1} \sum\limits_{(i_1,\cdots,
i_n)\in\Theta_{n}^{t}}\lambda ^{t-1}D^{i_1}(v_1)\cdots D^{i_n}(v_n)
\end{equation}
Note that in (\ref{e3}), $D(v_1)\cdots D(v_n)>D^{i_1}(v_1)\cdots
D^{i_n}(v_n)$.

 We prove the result by induction on $i$. For $i=0$ or 1, the result is clear.

 Now we assume that the result is true for $i-1,
i\geq 2$.   Since
$$
D^{i}(v_1v_2\cdots v_n)=D(D^{i-1}(v_1v_2\cdots v_n)),
$$
by (\ref{e3}) and by induction, we have
$$
\overline{D^{i}(v_1v_2\cdots v_n)}=D^{i}(v_1)D^{i}(v_2)\cdots
D^{i}(v_n)
$$
and the coefficient of $\overline{D^{i}(v_1v_2\cdots v_n)}$ is
$\lambda^{(n-1)i}$.\hfill $\blacksquare$

\ \

Let $D\langle X;\Omega\rangle$ be a free $\lambda$-differential
associative algebra with multiple operators  $\Omega$ on $X$ and
$\star\notin X$. By a $\star$-$\Omega$-word we mean any expression
in $\mathfrak{L}(X\cup \{\star\})$ with only one occurrence of
$\star$. The set of all $\star$-$\Omega$-words on $X$ is denoted by
$\mathfrak{L}^\star (X)$.

Let $u$ be a $\star$-$\Omega$-word and $s\in D\langle
X;\Omega\rangle$. Then we call
$$
u|_{s}=u|_{\star\mapsto s}
$$
an $s$-$\Omega$-word.

In other words, an $s$-$\Omega$-word $u|_{\star\mapsto s}$ means
that we have replaced the $\star$ of $u$ by $s$.

For example, if
$u=\theta_3(D(x_1)D^2(x_2),D^2(\star),D^3(x_4)D^3(x_5))$, then

 \setlength{\unitlength}{1cm}
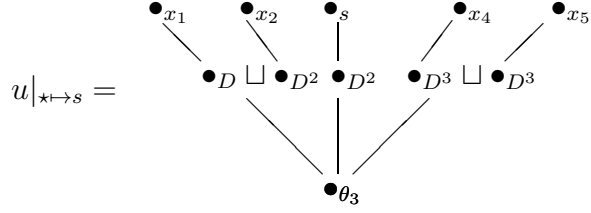
\begin{picture}(4,3.0)
\put(2.5,1.5){$u|_{\star\mapsto s} =$} \put(7,0.5){\line(1,1){1}}
\put(6.8,0.5){\line(0,1){1}} \put(6.6,0.5){\line(-1,1){1}}
\put(6.6,0.2){$\bullet_{\theta_3}$}
\put(6.6,0.2){$\bullet_{\theta_3}$}\put(6.7,1.7){$\bullet_{D^2}$}
\put(5,1.7){$\bullet_{D}\sqcup\bullet_{D^2}$}\put(7.7,1.7){$\bullet_{D^3}\sqcup\bullet_{D^3}$}
\put(5.0,2.0){\line(-1,1){0.5}}\put(4.3,2.6){$\bullet_{x_1}$}
\put(6,2.0){\line(-3,4){0.4}}\put(5.5,2.6){$\bullet_{x_2}$}\put(6.8,2.0){\line(0,1){0.5}}
\put(6.6,2.6){$\bullet_{s}$}\put(7.9,2.0){\line(3,4){0.4}}\put(8.3,2.6){$\bullet_{x_4}$}
\put(9.0,2.0){\line(1,1){0.5}}\put(9.6,2.6){$\bullet_{x_5}$}
\end{picture}\\

Similarly, we can define  $(\star_1, \star_2)$-$\Omega$-words as
expressions in  $\mathfrak{L}(X\cup \{\star_1, \star_2\})$  with
only one occurrence of $\star_1$ and only one occurrence of
$\star_2$. Let us denote by $\mathfrak{L}^{\star_1, \star_2} (X)$
the set of all $(\star_1, \star_2)$-$\Omega$-words. Let $u\in
\mathfrak{L}^{\star_1, \star_2} (X)$. Then we call
$$
u|_{s_1, s_2}= u|_{\star_1\mapsto s_1,\star_2\mapsto s_2}
$$
an $s_1$-$s_2$-$\Omega$-word.

If $\overline{u|_{s}}=u|_{\overline{s}}$, then we call the
$s$-$\Omega$-word $u|_{s}$  a normal $s$-$\Omega$-word.

Note that not each $s$-$\Omega$-word is a normal $s$-$\Omega$-word,
for example,
 if
$$
 u=\theta_3(D(x_1) D^2(\star), D^2(x_3),D^3(x_4)D^3(x_5))\in
\mathfrak{L}^\star(X)
$$
and $s=xy+1$, then $u|_{s}$ is not a normal $s$-$\Omega$-word.
However, if we take
$$u'=\theta_3(D(x_1)\star, D^2(x_3),D^3(x_4)D^3(x_5))\in \mathfrak{L}^\star(X),$$
then $u|_{s}=u'|_{\star\ \mapsto D^2(s)}$ and $u'|_{\star\  \mapsto
D^2(s)}$ is a normal $D^2(s)$-$\Omega$-word. By this way, we can
easily prove the following lemma.

\begin{lemma}\label{l3.4}
For any $s$-word $u|_{s}$, there exist $i\geq 0$ and $u'\in
\mathfrak{L}^\star(X)$ such that $u|_{s}=u'|_{D^{i}(s)}$ and
$u'|_{D^{i}(s)}$ is  a normal $D^{i}(s)$-word.\hfill $\blacksquare$
\end{lemma}

By Lemmas \ref{l3.1},  \ref{l3.3}, we have

\begin{lemma}\label{l2}
For any $u, v \in \mathfrak{L}(X), w \in \mathfrak{L}^\star(X), \ \
u>v\Rightarrow \overline{w|_u}>\overline{w|_v}. $\hfill
$\blacksquare$
\end{lemma}

 Let $f, g \in D\langle X;\Omega\rangle$ be monic. Then,
there are two kinds of compositions.

\begin{enumerate}
\item[(i)]If there exists a $w=\overline{D^{i}(f)} a=b \overline{D^{j}(g)}$ for some $a,b\in
\mathfrak{L}(X)$ such that $bre(w)< bre(\bar{f})+bre(\bar{g})$, then
we call
$(f,g)_{w}=lc(D^{i}(f))^{-1}D^{i}(f)a-lc(D^{j}(g))^{-1}bD^{j}(g)$
the intersection composition of $f$ and $g$ with respect to $w$.

\item[(ii)] If there exists a $u \in \mathfrak{L}^\star(X)$ such that
$w=\overline{D^{i}(f)}=u|_{\overline{D^{^{j}}({g})}}$, then we call
$(f,g)_{w}=lc(D^{i}(f))^{-1}D^{i}(f)-lc(D^{j}(g))^{-1}u|_{D^{^{j}}({g})}$
the including composition of $f$ and $g$ with respect to $w$.
\end{enumerate}

\noindent {\bf Remark:} In (i) and (ii),
$lc(D^{i}(f))=lc(D^{j}(g))=1$ if $\lambda=0$;
$lc(D^{i}(f))=\lambda^{(bre(\bar{f})-1)i}$ and
$lc(D^{j}(g))=\lambda^{(bre(\bar{g})-1)j}$  if $\lambda\neq0$.

In the above definition, $w$ is called the ambiguity of the
composition. Clearly,
$$
(f,g)_w\in Id(f,g)
$$
where $Id(f,g)$ is the ideal of $D\langle X;\Omega\rangle$ generated
by $f,\ g$.

It is noted that by an ideal $I$ of $D\langle X;\Omega\rangle$ we
mean $I$ is an ideal of $D\langle X;\Omega\rangle$ as associative
algebra and closed under $D$ and $\Omega$.

By Lemma \ref{l2}, we have
$$
\overline{(f,g)_w}< w.
$$

Let $f, g \in D\langle X;\Omega\rangle$  with
$\bar{f}=u|_{\overline{D^{i}(g)}}$ for some  $u\in
\mathfrak{L}^*(X)$. Then the transformation
$$
f\rightarrow f-lc(f)  u|_{lc(D^{i}(g))^{-1}D^{i}(g)}
$$
is called the elimination of the leading word (ELW) of $f$ by $g$.

Let $S$ be a monic subset of  $D\langle X;\Omega\rangle$. Then the
composition $(f,g)_w$ is called trivial modulo $(S,w)$ if
$$
(f,g)_w=\sum\alpha_iu_i|_{D^{l_i}(s_i)}
$$
where each $\alpha_i\in K$,  $u_i\in \mathfrak{L}^\star(X)$, $s_i\in
S$, $u_i|_{D^{l_i}(s_i)}$ is a normal $D^{l_i}(s_i)$-word and
$u_i|_{\overline{D^{l_i}(s_i)}}< w$. If this is the case, we write
$$
(f,g)_w\equiv 0 \ \ mod (S,w).
$$
In general, for any two polynomials $p$ and $q$, $ p\equiv q \ \ mod
(S,w) $ means that $ p-q=\sum\alpha_iu_i|_{D^{l_i}(s_i)}, $ where
each $\alpha_i\in K$,  $u_i\in \mathfrak{L}^\star(X)$, $s_i\in S$,
$u_i|_{D^{l_i}(s_i)}$ is a  normal $D^{l_i}(s_i)$-word and
$u_i|_{\overline{D^{l_i}(s_i)}}< w$.

$S$ is called a Gr\"{o}bner-Shirshov basis  in  $D\langle
X;\Omega\rangle$ if any composition $(f,g)_w$ of $f,g\in S$  is
trivial modulo $(S,w)$.

\begin{lemma}\label{l3.5}
Let $S$ be a Gr\"{o}bner-Shirshov  basis  in  $D\langle
X;\Omega\rangle, \ u_1, \ u_2\in \mathfrak{L}^\star(X)$
 and $s_1, s_2\in S$. If
 $w=u_1|_{_{\overline{D^{l_1}(s_1)}}}=u_2|_{_{\overline{D^{l_2}(s_2)}}}$,
 where  each $u_i|_{_{D^{l_i}(s_i)}}$ is  a  normal $D^{l_i}(s_i)$-word, $i=1,2$,  then
$$
u_1|_{lc(D^{l_1}(s_1))^{-1}D^{l_1}(s_1)}\equiv
u_2|_{lc(D^{l_2}(s_2))^{-1}D^{l_2}(s_2)} \ \ mod (S,w).
$$
\end{lemma}
{\bf Proof:} \ There are three cases to consider.

(i)\ \  $\overline{D^{l_1}(s_1)}$ and $\overline{D^{l_2}(s_2)}$ are
disjoint.  Then there  exits a $(\star_1,\star_2)$-$\Omega$-word
$\Pi$ such that
$$\Pi|_{_{\overline{D^{l_1}(s_1)},\
\overline{D^{l_2}(s_2)}}}=u_1|_{_{\overline{D^{l_1}(s_1)}}}=u_2|_{_{\overline{D^{l_2}(s_2)}}}.
$$
Then
\begin{eqnarray*} && u_2|_{_{lc(D^{l_2}( s_2))^{-1}D^{l_2}( s_2)}}-u_1|_{_{lc(D^{l_1}
(s_1))^{-1}D^{l_1}
(s_1)}}\\
&=&\Pi|_{\overline{D^{l_1}(s_1)}, \  lc(D^{l_2}( s_2))^{-1}D^{l_2} (s_2)}-
\Pi|_{lc(D^{l_1}(s_1))^{-1}D^{l_1}(s_1), \ \overline{D^{l_2}(s_2)}}\\
&=&-\Pi|_{lc(D^{l_1}
(s_1))^{-1}D^{l_1}(s_1)-\overline{D^{l_1}(s_1)}, \ lc(D^{l_2}(
s_2))^{-1} D^{l_2} (s_2)}\\
&&+\Pi|_{lc(D^{l_1} (s_1))^{-1}D^{l_1}(s_1),\ lc(D^{l_2} (s_2))^{-1}
D^{l_2} (s_2)-\overline{D^{l_2} (s_2)}}.
\end{eqnarray*}

 Let
\begin{eqnarray*}
-\Pi|_{lc(D^{l_1} (s_1))^{-1}D^{l_1}(s_1)-\overline{D^{l_1}(s_1)}, \
lc(D^{l_2}( s_2))^{-1} D^{l_2}
(s_2)}&=&\sum\alpha_{2_t}u_{2_t}|_{_{D^{l_{2}}(s_2)}},
\\
\Pi|_{lc(D^{l_1} (s_1))^{-1}D^{l_1}(s_1),\ lc(D^{l_2} (s_2))^{-1}
D^{l_2} (s_2)-\overline{D^{l_2} (s_2)}}
&=&\sum\alpha_{1_l}u_{1_l}|_{_{D^{l_1}(s_1)}}
\end{eqnarray*}
where all $u_{2_t}|_{D^{l_{2}}(s_2)}$ and $u_{1_l}|_{D^{l_1}(s_1)}$
are normal $D^{l_{2}}(s_2)$- and $D^{l_1}(s_1)$-words, respectively.

By Lemma \ref{l2}, we have
$\overline{lc(D^{l_i}(s_i)^{-1})D^{l_i}(s_i)-\overline{D^{l_i}(s_i)}}<\overline{D^{l_i}(s_i)},
\ i=1,2$.   It follows that
$$
u_1|_{lc(D^{l_1}(s_1))^{-1}D^{l_1}(s_1)}\equiv
u_2|_{lc(D^{l_2}(s_2))^{-1}D^{l_2}(s_2)} \ \ mod (S,w).
$$

(ii)   $\overline{D^{l_1}(s_1)}$ and $\overline{D^{l_2}(s_2)}$ have
nonempty intersection but do not
 include each other. Without lost of generality we can assume that
$\overline{D^{l_1}(s_1)}a=b\overline{D^{l_2}(s_2)}$ for some $a,
b\in\mathfrak{L}(X)$. Then there exists a
$\Pi\in\mathfrak{L}^\star(X)$ such that
$$
\Pi|_{\overline{D^{l_1}(s_1)}a}=u_1|_{\overline{D^{l_1}(s_{_1}})}=
u_2|_{\overline{D^{l_2}(s_2)}}=\Pi|_{b\overline{D^{l_2}(s_2)}}
$$
where $\Pi|_{D^{l_1}(s_1)a}$ is a  normal $D^{l_1}(s_1)a$-word and
$\Pi|_{bD^{l_2}(s_2)}$ is a normal $bD^{l_2}(s_2)$-word.   Thus, we
have
\begin{eqnarray*}
&&u_{_2}|_{lc(D^{l_2}( s_{_2}))^{-1}D^{l_2}( s_{_2})}-u_1|_{
lc(D^{l_1}( s_{_1}))^{-1}D^{l_1}(s_{_1})}\\
&=&\Pi|_{lc(D^{l_2}( s_{_2}))^{-1}bD^{l_2}(
s_{_2})}-\Pi|_{lc(D^{l_1}(
s_{_1}))^{-1}D^{l_1}(s_{_1})a}\\
&=&-\Pi|_{lc(D^{l_1}( s_{_1}))^{-1}D^{l_1}(s_{_1})a-lc(D^{l_2}(
s_{_2}))^{-1}bD^{l_2}(s_{_2})}.
\end{eqnarray*}
Since  $S $ is a Gr\"{o}bner-Shirshov basis in  $D\langle
X;\Omega\rangle$,
 we have
$$
lc(D^{l_1}( s_{_1}))^{-1}D^{l_1}(s_{_1})a-lc(D^{l_2}(
s_{_2}))^{-1}bD^{l_2}(s_{_2})=\sum\alpha_jv_j|_{D^{t_j}(s_j)}
$$
where each $\alpha_j\in K, \ v_j\in \mathfrak{L}^\star(X), \ s_j\in
S$, $v_j|_{\overline{D^{t_j}(s_j)}}<\overline{ D^{l_1}(s_1)}a$ and
$v_j|_{D^{t_j}(s_j)}$ is a normal  $D^{t_j}(s_j)$-word.  Let $
\Pi|_{v_j|_{D^{t_j}(s_j)}}=\Pi_{j}|_{D^{t_j}(s_j)} $. Then
$\Pi_{j}|_{D^{t_j}(s_j)}$ is also a normal $D^{t_j}(s_j)$-word.
Therefore
$$
u_{_2}|_{lc(D^{l_2}( s_{_2}))^{-1}D^{l_2}( s_{_2})}-u_1|_{
lc(D^{l_1}(
s_{_1}))^{-1}D^{l_1}(s_{_1})}=\sum\alpha_j\Pi_{j}|_{D^{t_j}(s_j)}
$$
with $ \Pi_{j}|_{\overline{D^{t_j}(s_j)}}<w. $ It follows that
$$
u_1|_{lc(D^{l_1}(s_1))^{-1}D^{l_1}(s_1)}\equiv
u_2|_{lc(D^{l_2}(s_2))^{-1}D^{l_2}(s_2)} \ \ mod (S,w).
$$

(iii) One of  $\overline{D^{l_1}(s_1)}$, $\overline{D^{l_2}(s_2)}$
is contained in the other. For example, let $
\overline{D^{l_1}(s_1)} =u|_{\overline{D^{l_2}(s_2)}}$ for some
$\star$-word $u$. Then
$$
w=u_2|_{\overline{D^{l_2}(s_2)}}=u_1|_{u|_{\overline{D^{l_2}(s_2)}}}\
\ \mbox{ and }
$$
\begin{eqnarray*}
&&u_2|_{
lc(D^{l_2}(s_2))^{-1}D^{l_2}(s_2)}-u_1|_{lc(D^{l_1}(s_1))^{-1}
D^{l_1}(s_1)}\\
&=&u_1|_{ lc(D^{l_2}(s_2))^{-1}u|_{D^{l_2}(s_2)}}-u_1|_{lc(D^{l_1}(s_1))^{-1}D^{l_1}(s_1)}\\
&=&-u_1|_{ lc(D^{l_1}(s_1))^{-1}D^{l_1}(s_1)-lc(D^{l_2}(s_2))^{-1}
u|_{D^{l_2}(s_2)}}.
\end{eqnarray*}
Similar to  (ii), we can obtain the result. \hfill $\blacksquare$

The following theorem is an analogy  of Shirshov's Composition lemma
for Lie algebras \cite{S3}, which was specialized to associative
algebras by Bokut \cite{b76}, see also Bergman \cite{b}.

\begin{theorem}\label{3.6}{\em(Composition-Diamond lemma)}\ \  Let $S$ be
a monic  subset of  $D\langle X;\Omega\rangle$, $Id(S)$ the ideal of
 $D\langle X;\Omega\rangle$ generated by $S$ and
 $>$ the order on $\mathfrak{L}(X)$ defined as before.  Then the following
statements are equivalent:
 \begin{enumerate}
\item[(I)] $S $ is a Gr\"{o}bner-Shirshov basis in $D\langle X;\Omega\rangle$.
\item[(II)] $ f\in Id(S)\Rightarrow
\bar{f}=u|_{\overline{D^{i}(s)}}$  for some $u \in
\mathfrak{L}^\star(X)$, $s\in S$ and $i\geq 0$.
 \item[(III)] $Irr(S) = \{ w\in \mathfrak{L}(X) |  w \neq
u|_{\overline{D^{i}(s)}},\ s\in S,\ u\in\mathfrak{L}^\star(X), \
u|_{D^{i}(s)}
    \mbox{ is a normal  } D^{i}(s)\mbox{-word}\}$
is a $K$-basis of $D\langle X;\Omega|S\rangle=D\langle
X;\Omega\rangle/Id(S)$.
\end{enumerate}
\end{theorem}
{\bf Proof:}\ \
 (I)$\Longrightarrow$ (II)\ \ Let  $0\neq f\in
Id(S)$. Then by Lemma \ref{l3.4}, we can assume that
$$
f=\sum\limits_{i=1}^{n}\alpha_i u_i|_{_{lc(D^{l_i}
(s_i))^{-1}D^{l_i} (s_i)}}
$$
 where each $\alpha_i\in K$,  $u_i\in \mathfrak{L}^\star(X)$, $s_i\in
S$ and $u_i|_{D^{l_i}(s_i)}$ is a normal $D^{l_i}(s_i)$-word. Let
$w_i=u_i|_{\overline{D^{l_i}( s_i)}}$. We arrange these leading
terms in non-increasing order by
$$
w_1= w_2=\cdots=w_m >w_{m+1}\geq \cdots\geq w_n.
$$
Now we prove the result by  induction  on $m$.

If $m=1$, then $\bar{f}=u_1|_{\overline{D^{l_1} (s_1)}}$.

Now we assume that $m\geq 2$. Then $ u_1|_{\overline{D^{l_1}(
s_1)}}=w_1=u_2|_{\overline{D^{l_2}( s_2)}}$.

We prove the result by induction on $w_1$.  Since $S $ is a
Gr\"{o}bner-Shirshov basis  in $D\langle X;\Omega\rangle$, by Lemma
\ref{l3.5}, we have
$$
u_2|_{lc( D^{l_2}(s_2))^{-1} D^{l_2}(s_2)}-u_1|_{lc(
D^{l_1}(s_1))^{-1} D^{l_1}( s_1)}=\sum\beta_jv_j|_{D^{t_j}(s_j)}
$$
where $\beta_j\in K, \ s_j\in S, v_j\in \mathfrak{L}^\star(X)$,
$v_j|_{\overline{D^{t_j}(s_j)}}<w_1$ and $v_j|_{D^{t_j}(s_j)}$ is a
normal $D^{t_j}(s_j)$-word. Therefore, since
\begin{eqnarray*}
&& \alpha_1u_1|_{lc(D^{l_1}( s_1))^{-1}D^{l_1}(
s_1)}+\alpha_2u_2|_{lc(D^{l_2}( s_2))^{-1} D^{l_2}(
s_2)}\\
&=&(\alpha_1+\alpha_2)u_1|_{lc(D^{l_1}( s_1))^{-1}D^{l_1}(
s_1)}+\alpha_2(u_2|_{lc(D^{l_2}( s_2))^{-1} D^{l_2}(
s_2)}-u_1|_{lc(D^{l_1}( s_1))^{-1}D^{l_1}( s_1)}),
\end{eqnarray*}
we have
\begin{eqnarray*}
f=(\alpha_1+\alpha_2)u_1|_{lc(D^{l_1}( s_1))^{-1} D^{l_1}(
s_1)}+\sum\alpha_2\beta_jv_j|_{D^{t_j}(s_j)}+
\sum\limits_{i\geq3}\alpha_iu_i|_{lc(D^{l_i}(
s_i))^{-1}D^{l_i}(s_i)}.
\end{eqnarray*}

If either $m>2$ or $\alpha_1+\alpha_2\neq 0$, then the result
follows from induction on $m$. If $m=2$ and $\alpha_1+\alpha_2=0$,
then  the result follows from induction on $w_1$.

(II)$\Longrightarrow$ (III) For any $f\in D\langle X;\Omega\rangle$,
by induction on $\bar{f}$,  we have
$$
f=\sum \limits_{u_i\in Irr(S), \ u_i\leq\bar{f}}\alpha_iu_i + \sum
\limits_{s_j\in S,\
v_j|_{_{\overline{D^{l_j}(s_j)}}}\leq\bar{f}}\beta_jv_j|_{_{D^{l_j}(s_j)}}
$$
where $\alpha_i,\beta_j\in K, \ v_j|_{D^{l_j}(s_j)}$ is a normal
$D^{l_j}(s_j)$-word.  Then $f+Id(S)$ can be expressed by the
elements of $Irr(S)$.
 Now suppose that $\alpha_1u_1+\alpha_2u_2+\cdots\alpha_nu_n=0$ in
$ D\langle X;\Omega|S\rangle$ with each $0\neq\alpha_i\in K,\ u_i\in
Irr(S)$, Then, in $ D\langle X;\Omega\rangle$,
$$
g=\alpha_1u_1+\alpha_2u_2+\cdots+\alpha_nu_n\in Id(S).
$$
By (II), we have $\bar{g}=u_i\notin Irr(S)$ for some $1\leq i\leq
n$, a contradiction. Hence, $Irr(S)$ is $K$-linearly independent.
This shows that $Irr(S)$ is a $K$-basis of $ D\langle
X;\Omega|S\rangle$.

(III)$\Longrightarrow $(I) For any composition $(f,g)_w$ in $S$,
since $(f,g)_w\in Id(S)$ and (III), we have
$$
(f,g)_w=\sum\beta_jv_j|_{_{D^{l_j}(s_j)}}
$$
where each $\beta_j\in K$,  $v_j\in \mathfrak{L}^\star(X)$, $s_j\in
S$ and $v_j|_{_{\overline{D^{l_j}(s_j)}}}\leq\overline{(f,g)_w}< w$.
This shows (I). \hfill $\blacksquare$

\section{ Gr\"{o}bner-Shirshov bases for free  $\lambda$-differential Rota-Baxter  algebras}

In this section, we obtain a Gr\"{o}bner-Shirshov basis and a linear
basis for a free $\lambda$-differential Rota-Baxter algebra.
Consequently,  we construct the  free  $\lambda$-differential
Rota-Baxter  algebra on  set $X$.

Let $K$ be a field and $\lambda\in K$.
 A differential Rota-Baxter algebra of weight $\lambda$
 (\cite{lw}), called also $\lambda$-differential Rota-Baxter
 algebra,
is an associative $K$-algebra $R$ with two $K$-linear operators
$P,D:R\rightarrow R$ such that for any $u,v\in R$,

(i) (Rota-Baxter relation) $P(u)P(v)=P(uP(v))+P(P(u)v)+\lambda
P(uv);$

(ii) ($\lambda$-differential relation) $ D(uv)=\lambda
D(u)D(v)+D(u)v+uD(v);$

(iii) $D(P(u))=u$.

Hence, any  $\lambda$-differential Rota-Baxter
 algebra is a    $\lambda$-differential  associative algebra with operator  $\Omega=\{P\}$.

Throughout this section, we assume $\Omega=\{P\}$ and $D\langle
X;\Omega \rangle$ the free $\lambda$-differential associative
algebra with operator $\Omega$.

 Let $S$ be the subset of $D\langle X;\Omega \rangle$ consisting of  the following
polynomials:
\begin{enumerate}
\item[1.]
$P(u)P(v)-P(uP(v))-P(P(u)v)-\lambda P(uv), \ u,\ v\in
\mathfrak{L}(X)$,
\item[2.] $D(P(u))-u, \ u\in \mathfrak{L}(X)$.
\end{enumerate}

 Denoted by
$$f(u,v)=P(u)P(v)-P(uP(v))-P(P(u)v)-\lambda P(uv), \ u,\ v\in \mathfrak{L}(X).$$
\begin{lemma} \label{l4.1} Let  $j> 0$.

(i) If   $\lambda\neq 0$, then  $mod(S, D^{j}(P(u))D^{j}(P(v)))$
\begin{eqnarray*}
D^{j}(f(u,v)) \equiv
\lambda^{j}(D^{j}(P(u))D^{j}(P(v))-D^{j-1}(u)D^{j-1}(v)).
\end{eqnarray*}

(ii) If   $\lambda=0$, then  $mod(S, D^{j}(P(u))P(v))$
\begin{eqnarray*}
D^{j}(f(u,v)) \equiv D^{j}(P(u))P(v)-D^{j-1}(u)P(v).
\end{eqnarray*}
\end{lemma}
{\bf Proof: }We just prove (i). The proof of (ii) is similar to (i).

 Induction on $j$. For $j=1$, we have
\begin{eqnarray*}
&& D(f(u,v))\\
&=& D(P(u)P(v))-D(P(uP(v)))-D(P(P(u)v))-\lambda D(P(uv))\\
&=&\lambda
D(P(u))D(P(v))+D(P(u))P(v)+P(u)D(P(v))+\\
&&-D(P(uP(v)))
-D(P(P(u)v))-\lambda D(P(uv)) \\
&\equiv& \lambda
D(P(u))D(P(v))+uP(v)+P(u)v -uP(v)-P(u)v-\lambda uv\\
&\equiv& \lambda (D(P(u))D(P(v))-uv)\ \ mod(S,D(P(u))D(P(v)).
\end{eqnarray*}
Now we assume that the result is true for $j-1, j\geq 1$, i.e.,
$$
D^{j-1}(f(u,v))=\lambda^{j-1}
(D^{j-1}(P(u))D^{j-1}(P(v))-D^{j-2}(u)D^{j-2}(v))+\sum \alpha_i
u_i|_{D^{t_i}(s_i)}
$$
where each $\alpha_i\in K$, $s_i\in
S$,$u_i|_{\overline{D^{t_i}(s_i)}}<D^{j-1}(P(u))D^{j-1}(P(v))$ and
$u_i|_{D^{t_i}(s_i)}$ is a  normal $D^{t_i}(s_i)$-word. Therefore
$$
D(\sum \alpha_i u_i|_{D^{t_i}(s_i)})=\sum
\beta_{l}v_l|_{D^{k_l}(s_l)}
$$
where each $v_l|_{D^{k_l}(s_l)}$ is normal $D^{k_l}(s_l)$-word and
$v_l|_{ \overline{D^{k_l}(s_l) }}<
\overline{D(D^{j-1}(P(u))D^{j-1}(P(v)))}=D^{j}(P(u))D^{j}(P(v))$  by
Lemma \ref{l3.3}. So, $mod(S, D^{j}(P(u))D^{j}(P(v)))$, we have
\begin{eqnarray*}
&&D^{j}(f(u,v))=D(D^{j-1}(f(u,v)))\\
&\equiv& \lambda^{j-1}D(
D^{j-1}(P(u))D^{j-1}(P(v))-D^{j-2}(u)D^{j-2}(v))\\
&\equiv& \lambda^{j-1}(\lambda D^{j}(P(u))D^{j}(P(v))+D^{j}(P(u))D^{j-1}(P(v))+D^{j-1}(P(u))D^{j}(P(v)))\\
&&- \lambda^{j-1}(\lambda D^{j-1}(u)D^{j-1}(v)+D^{j-1}(u)D^{j-2}(v)+D^{j-2}(u)D^{j-1}(v) )\\
&\equiv& \lambda^{j-1}(\lambda D^{j}(P(u))D^{j}(P(v))+ D^{j-1}(u)D^{j-2}(v)+D^{j-2}(u)D^{j-1}(v))\\
&&- \lambda^{j-1}(\lambda D^{j-1}(u)D^{j-1}(v)+D^{j-1}(u)D^{j-2}(v)+D^{j-2}(u)D^{j-1}(v) )\\
&\equiv& \lambda^{j}(D^{j}(P(u))D^{j}(P(v))-D^{j-1}(u)D^{j-1}(v)). \
\ \ \ \ \ \ \ \ \ \ \ \ \ \ \ \ \ \ \ \ \ \ \ \ \ \ \ \ \ \ \ \ \ \
\ \ \ \ \ \ \ \ \  \ \ \  \ \hfill \blacksquare
\end{eqnarray*}

\begin{theorem}\label{t4.1}With  the order $>$
 on $\mathfrak{L}(X)$ defined as before,  $S$ is a Gr\"{o}bner-Shirshov
basis in $D\langle X;\Omega\rangle$.
\end{theorem}
\noindent {\bf Proof.} Denote by $k\wedge l$ the composition of the
polynomials of type $k$ and type $l$. There are two cases $\lambda
\neq 0$ and $\lambda =0$ to consider.

(i)  For $\lambda \neq 0$, the ambiguities of all possible
compositions of the polynomials  in $S$
 are list as below:
 \begin{tabbing}
  $2\wedge2\ \ \ D^{i}(D(P(u|_{D^{^{j}}(D(P(v)))}))),$\\[0.7ex]
  $2\wedge1\ \ \  D^{i}(D(P(u|_{D^{^{j}}(P(v))D^{^{j}}(P(w))}))),$\\[0.7ex]
   $1\wedge2\ \ \ D^{j}(P(u|_{D^{^{i}}(D(P(v)))}))D^{j}(P(w)), \
   D^{j}(P(v))D^{j}(P(u|_{D^{^{i}}(D(P(w)))})),$\\[0.7ex]
 $1\wedge1\ \ \ D^{^{j}}(P(u))D^{^{j}}(P(v))D^{^{j}}(P(w)),\
 D^{j}(P(u|_{D^{i}(P(v))D^{i}(P(w))}))D^{j}(P(v')),$\\[0.7ex]
 $\ \ \ \ \ \ \  \ \ D^{j}(P(v))D^{j}(P(u|_{D^{i}(P(w))D^{i}(P(v'))}))$
\end{tabbing}
where $u\in\mathfrak{L}^\star(X) ,v,w, v'\in \mathfrak{L}(X)$ and
$i,j\geq 0$.
 Now we check that all the compositions are trivial.
\begin{eqnarray*}
2\wedge2&=& D^{i}(D(P(u|_{D^{^{j}}(D(P(v)))}))-u|_{D^{^{j}}(D(P(v)))})-D^{i}(D(P(u|_{D^{^{j}}(D(P(v)-v))})))\\
&=& -D^{i}(u|_{D^{^{j}}(D(P(v)))})+D^{i}(D(P(u|_{D^{^{j}}(v))})))\\
&\equiv& -D^{i}(u|_{D^{^{j}}(v)})+D^{i}(u|_{D^{^{j}}(v)})\\
&\equiv& 0.
\end{eqnarray*}

Now, assume that $j>0$.
\begin{eqnarray*}
2\wedge1&=& D^{i}(D(P(u|_{D^{^{j}}(P(v))D^{^{j}}(P(w))}))) -
D^{i}(u|_{D^{^{j}}(P(v))D^{^{j}}(P(w))})\\
 &&-\lambda^{-j}D^{i}(D(P(u|_{D^{^{j}}(f(v,w))})))\\
&\equiv& D^{i}(D(P(u|_{D^{^{j}}(P(v))D^{^{j}}(P(w))})))-D^{i}(u|_{D^{^{j}}(P(v))D^{^{j}}(P(w))})\\
& &
-\lambda^{-j}D^{i}(D(P(u|_{\lambda^{^j}(D^{j}(P(v))D^{j}(P(w))-D^{j-1}(v)D^{j-1}(w))})))\
\
\mbox{(by Lemma \ref{l4.1})}\\
&\equiv& -D^{i}(u|_{D^{^{j}}(P(v))D^{^{j}}(P(w))})+ D^{i}(D(P(u|_{D^{j-1}(v)D^{j-1}(w))})) \\
&\equiv&
-D^{i}(u|_{D^{^{j-1}}(v)D^{^{j-1}}(w)})+D^{i}(u|_{D^{^{j-1}}(v)D^{^{j-1}}(w)})\\
&\equiv& 0.
\end{eqnarray*}
\begin{eqnarray*}
1\wedge2&=& \lambda^{-j}D^{j}(f(u|_{D^{^{i}}(D(P(v)))},w))-D^{j}(P(u|_{D^{^{i}}(D(P(v))-v)})D^{j}(P(w))\\
&\equiv&\lambda^{-j}
(\lambda^{j}(D^{j}(P(u|_{D^{^{i}}(D(P(v)))}))D^{j}(P(w))-D^{j-1}(u|_{D^{^{i}}(D(P(v)))})D^{j-1}(w)))\\
&&-D^{j}(P(u|_{D^{^{i}}(D(P(v)))})D^{j}(P(w))+D^{j}(P(u|_{D^{^{i}}(v)})D^{j}(P(w)) \ \
\mbox{(by Lemma \ref{l4.1})}\\
&\equiv& 0.
\end{eqnarray*}
Similarly, we can prove another case of $1\wedge2$ to be also
trivial.
\begin{eqnarray*}
1\wedge1&=& \lambda^{-j}D^{j}(f(u,v))D^{j}(P(w))-\lambda^{-j}D^{j}(P(u))D^{j}(f(u,w))\\
&\equiv&
(D^{j}(P(u))D^{j}(P(v))-D^{j-1}(u)D^{j-1}(v))D^{j}(P(w))\\
&&-D^{j}(P(u))(D^{j}(P(v))D^{j}(P(w))-D^{j-1}(v)D^{j-1}(w)) \
\mbox{(by Lemma \ref{l4.1})}\\
&\equiv&-D^{j-1}(u)D^{j-1}(v))D^{j-1}(w)+D^{j-1}(u)D^{j-1}(v))D^{j-1}(w)\\
&\equiv& 0.
\end{eqnarray*}

Another two cases in $1\wedge1$ can be easily checked.

For $j=0$, the proof is similar. We omit the details.

(ii) For $\lambda = 0$, the ambiguities of all possible compositions
of the polynomials  in $S$
 are list as below:
 \begin{tabbing}
  $2\wedge2\ \ \ D^{i}(D(P(u|_{D^{^{j}}(D(P(v)))}))),$\\[0.7ex]
  $2\wedge1\ \ \  D^{i}(D(P(u|_{D^{^{j}}(P(v))P(w)}))),$\\[0.7ex]
   $1\wedge2\ \ \ D^{j}(P(u|_{D^{^{i}}(D(P(v)))}))P(w), \
   D^{j}(P(v))P(u|_{D^{^{i}}(D(P(w))}),$\\[0.7ex]
 $1\wedge1\ \ \ D^{^{j}}(P(u))P(v)P(w),\
 D^{j}(P(u|_{D^{i}(P(v))P(w)}))P(v'), D^{j}(P(v))P(u|_{D^{i}(P(w))P(v')})$
\end{tabbing}
where $u\in\mathfrak{L}^\star(X) ,v,w, v'\in \mathfrak{L}(X)$ and
$i,j\geq 0$.

Similar to the proof of (i), we can prove that all compositions are
trivial modulo $S$. $\hfill \blacksquare$

\ \

 Let $Y$ and $Z$ be two subsets of
$\mathfrak{L}(X)$. Define the alternating product of $Y$ and $Z$
(see also \cite{kl1}):
$$
\Lambda(Y,Z)=(\cup_{ r\geq 1}(YP(Z))^r)\cup (\cup_{r\geq
0}(YP(Z))^rY) \cup(\cup_{ r\geq
1}(P(Z)Y)^r)\cup(\cup_{r\geq0}(P(Z)Y)^rP(Z)).
$$

Define
$$
\Phi_0=S(D^{\omega}(X)),
$$
and for $n>0$,
$$
\Phi_n=\Lambda(\Phi_0, \Phi_{n-1}).
$$
Let
$$
\Phi(D^{\omega}(X))=\bigcup_{n\geq 0 }\Phi_n.
$$

By   Theorem \ref{t4.1} and Theorem \ref{3.6}, we have the following
theorems.

\begin{theorem}{\em (\cite{b08})}
$Irr(S)= \Phi(D^{\omega}(X))$  is a $K$-basis of  $D\langle
X;\Omega|S\rangle$.
\end{theorem}

\begin{theorem}\label{t4.4}{\em(\cite{b08})}
$D\langle X;\Omega|S\rangle$ is a free $\lambda$-differential
Rota-Baxter algebra on  set $X$ with a linear basis
$\Phi(D^{\omega}(X))$.
\end{theorem}

We note here that Theorem \ref{t4.4} is also obtained by L. Guo and
W. Keigher in \cite{lw} by using the planar decorated rooted trees.

\ \

\noindent{\bf Acknowledgement}: The authors would like to thank
Professor L.A. Bokut for his guidance, useful discussions and
enthusiastic encouragement in writing up this paper.

\end{document}